\documentstyle[amstex,amssymb,A4,11pt]{article}
\makeatletter
\def\@@N{{\mathrm I \mkern-2.5mu \nonscript\mkern-.5mu \mathrm I
\mkern-5.5mu\mathrm  N}}
\def\@@Z{{
  \setbox0\hbox{\sf Z}
  \setbox1\hbox{\mathrm\kern.05\wd0
                \rlap{\vrule height.93\ht0 depth-.75\ht0 width.056\wd0 }%
                \kern-.13\wd0 \copy0 \kern-.6\wd0 \copy0 \kern-.1\wd0
                \llap{\vrule height.25\ht0 depth\z@ width.056\wd0}%
                \kern.05\wd0}
  \mathchoice{\copy1}{\copy1}{\mit Z\mkern-8mu Z}{\mit Z\mkern-7.5mu Z} }}
\def\@@insvline#1#2{{\setbox0\hbox{\m@th$#1\mathrm I$}
  \rlap{\m@th$#1 \mkern 5mu
  \vrule height.92\ht0 depth-.05\ht0 width.09\ht0 $}
  {\mathrm #2} }}
\def\@@Q{\mathpalette\@@insvline{Q}}
\def\@@R{{\mathrm I \mkern-2.5mu \nonscript\mkern-.5mu \mathrm R}}
\def\@@C{\mathpalette\@@insvline{C}}
\def\@@K{{\mathrm I \mkern-2.5mu \nonscript\mkern-.5mu\mathrm  K}}
\def\@@D{{\mathrm I \mkern-2.5mu \nonscript\mkern-.5mu\mathrm D}}

\def\callig#1{
  \ifcat#1a
    \ifnum`#1=\uccode`#1 {\cal #1}
    \else                #1
    \fi
  \else                  #1
  \fi}

\def\varbbbold#1{
  \ifcat#1a
    \ifnum`#1=\uccode`#1 \csname @@#1\endcsname
    \else                #1
    \fi
  \else                  #1
  \fi}

\def\eulerfraktur#1{
  \ifcat#1a {\frak #1}
  \else                  #1
  \fi}

\newif\if@@excloccurred
\newif\if@@questoccurred

\let\@@quest=? \catcode`?=\active
\let\@@excl=!  \catcode`!=\active \let\fact\@@excl

\def!{\protect\p@@doexcl}
\def?{\protect\p@@doquest}

\def\p@@doexcl{\let\@@xeq\relax\ifmmode\def\@@xeq{\futurelet\@@next\@@doexcl}
                               \else\@@excl\fi\@@xeq}
\def\p@@doquest{\let\@@xeq\relax\ifmmode\def\@@xeq{\futurelet\@@next\@@doquest}
                                \else\@@quest\fi\@@xeq}
\def\@@doexcl{\let\@@x@q\relax
     \if@@excloccurred \@@excloccurredfalse
                       \ifx\@@next\@sptoken\@@excl\@@excl
                       \else\let\@@x@q\doexclexcl\fi
\else\if@@questoccurred\@@questoccurredfalse
                       \ifx\@@next\@sptoken\@@quest\@@excl
                       \else\let\@@x@q\doquestexcl\fi
\else\ifx\@@next\@sptoken\@@excl
\else\ifx\@@next!\@@excloccurredtrue
\else\ifx\@@next?\@@excloccurredtrue
\else\let\@@x@q\doexcl\fi\fi\fi\fi\fi\@@x@q}

\def\@@doquest{\let\@@x@q\relax
     \if@@excloccurred \@@excloccurredfalse
                       \ifx\@@next\@sptoken\@@excl\@@quest
                       \else\let\@@x@q\doexclquest\fi
\else\if@@questoccurred\@@questoccurredfalse
                       \ifx\@@next\@sptoken\@@quest\@@quest
                       \else\let\@@x@q\doquestquest\fi
\else\ifx\@@next\@sptoken\@@quest
\else\ifx\@@next!\@@questoccurredtrue
\else\ifx\@@next?\@@questoccurredtrue
\else\let\@@x@q\doquest\fi\fi\fi\fi\fi\@@x@q}

\def\doexcl{\bold}        \def\doquest{\callig}
\def\doexclexcl{\varbbbold} \def\doexclquest{\boldsymbol}
\def\doquestexcl{\Bbb}    \def\doquestquest{\eulerfraktur}
\makeatother
\newcommand {\art}[6]{{\sc #1:} {#2.} {\em #3} {\bf #4}{(#5),} {#6.}}
\newcommand {\bo}[5]{{\sc #1:} {``#2."} {#3,} {#4} {#5.}}
\newcommand {\samp}[8]{{\sc #1:} {#2,} {{\em in} ``#3,"} {(#4, Ed.),}
                                  {pp. #5,} {#6,} {#7} {#8.}}

\newcommand {\unfart}[3]{{\sc #1:} {#2,} {(#3).}}
\newcommand {\prep}[1]{{\sc #1,}{\ in preparation.}}

\newcommand \force {{\hspace{0.4mm}}{\rule{0.1mm}{2.4mm}}
                       {\hspace{0.5mm}}{\rule{0.1mm}{2.4mm}}
               {\rule [1.2mm]{2.3mm}{0.1mm}}{\hspace{0.4mm}}}
\newcommand \reals {[\omega ]^\omega }
\newcommand \rval {]\!]\hspace{0.2mm}}
\newcommand \lval {[\![\hspace{0.2mm}}
\newcommand \one {{1\hspace{-0.4ex}{\rule{0.1ex}{1.5ex}}
             \hspace{-0.2ex}{\rule[1.5ex]{0.2ex}{0.1ex}}
             \hspace{-0.1ex}{\rule[1.6ex]{0.1ex}{0.05ex}}\hspace{0.5ex}}}

\newcommand \M {\operatorname{\mathbf{M}}}
\newcommand \pp {X}
\newcommand \DD {\mathfrak{D}}

\newcommand \Pd {{\operatorname{\mathbf{P}}}_D}
\newcommand \Pad {{\operatorname{\mathbf{P}}}_{\mathfrak{D}}}
\newcommand \Pf {{\operatorname{\mathbf{P}}}_{\mathfrak{F}}}
\newcommand \Pcf {{\operatorname{\mathbf{P}}}_{\cal{F}}}
\newcommand \Padname {{\operatorname{\mathbf{P}}}_{\tilde{\mathfrak{D}}}}

\newcommand \Pp {\operatorname{\mathbf{P}}}

\newcommand \Ww {\boldsymbol{\mathfrak{m}}}
\newcommand \fU {\boldsymbol{\mathfrak{U}}}
\newcommand \fF {\boldsymbol{\mathfrak{F}}}
\newcommand \fT {\boldsymbol{\mathfrak{T}}}
\newcommand \fw {\mathfrak{m}}

\newcommand \bU {\operatorname{\mathbf{U}}}
\newcommand \J  {{\cal{J}}}
\newcommand \LL {L}

\newcommand \ts {\tilde{s}}
\newcommand \tit {\tilde{t}}
\newcommand \tY {\tilde{Y}}
\newcommand \tX {\tilde{X}}
\newcommand \tsig {\tilde{\sigma}}
\newcommand \chs {\check{s}}
\newcommand \cht {\check{t}}
\newcommand \chX {\check{X}}
\newcommand \chY {\check{Y}}

\newcommand \nameDD {\tilde{\mathfrak{D}}}

\newcommand \oV {\omega_1^V}

\newcommand \ceq {\sqsubseteq}
\newcommand \kap {\sqcap}
\newcommand \puni {\sqcap}
\newcommand {\open}[2]{{(#1,} {#2)^{\omega}}}

\newcommand {\opf}[2]{{(#1,} {#2)^{\omega}_{\fF}}}
\newcommand {\join}[1]{{\stackrel{\frown}{\ }}{#1}}
\newcommand \parto {(\omega)^\omega}
\newcommand \partX {(X)^\omega}

\newcommand \NN {(!!N)}
\newcommand \mdom {\operatorname{dom}}
\newcommand \mmin {\operatorname{min}}
\newcommand \mmax {\operatorname{max}}
\newcommand \mMin {\operatorname{Min}}
\newcommand \cov {\operatorname{\mathbf{cov}}}
\newcommand \add {\operatorname{\mathbf{add}}}

\newcommand \ass {\stackrel{*}{=}}
\newcommand \seg {\preceq}
\newcommand \cp {\operatorname{cp}}
\newcommand \pac {\operatorname{pc}}
\newcommand \dR {{\mathfrak{R}}}
\newcommand \fD {{\mathfrak{D}}}
\newcommand \fG {{\mathfrak{G}}}
\newcommand \fJ {{\mathfrak{J}}}
\newcommand \fH {{\mathfrak{H}}}

\newcommand \dM {\operatorname{\boldsymbol{{\mathfrak{M}}}}}

\newcommand \Seq {\operatorname{Seq}}
\newcommand \la {\langle}
\newcommand \ra {\rangle}
\newcommand \GF {{\cal G}(\fF)}
\newcommand \cG {{\cal {G}}}
\newcommand \cO {{\cal O}}

\newcommand \subs {\subseteq}
\newcommand \cfF {\cap\fF}

\newtheorem {nummer}{ }[section]
\newtheorem {thm}[nummer]{T\footnotesize{\bf HEOREM}\bf }
\newtheorem {lm}[nummer]{L\footnotesize{\bf EMMA}\bf }

\newtheorem {fct}[nummer]{F\footnotesize{\bf ACT}\bf }

\newtheorem {cor}[nummer]{C\footnotesize{\bf OROLLARY}\bf }
\newcommand \rmk {{\bf{R{\footnotesize{\bf EMARK}}:\hspace*{3mm}}}}

\newcommand \proof {{\bf{P{\footnotesize{\bf ROOF}}:\hspace*{3mm}}}}

  \def\eop{{\unskip\nobreak\hfil\penalty50\hskip8mm\hbox{}
  \nobreak\hfil
  {$\boldsymbol{\dashv}$}\parfillskip=0mm \par\smallskip}}

\newcommand \eor {${\ }_\dashv$\hfill\smallskip}

\begin{document}

\begin{center}
      {\Large{\sc{Symmetries between two Ramsey properties}}}
\end{center}
\smallskip

\begin{center}
\small{Lorenz Halbeisen\footnote{The author
wishes to thank the {\em Swiss National Science Foundation\/} for
supporting him.}\\
        Universit\'e de Caen\\
        France}
\end{center}

\bigskip

\begin{abstract}
In this article we compare the
well-known Ramsey property with a dual form of
it, the so called dual-Ramsey property (which was suggested first
by Carlson and Simpson). Even if the two properties are different,
it can be shown that all classical results known for the Ramsey
property also hold for the dual-Ramsey property. We will also show
that the dual-Ramsey property is closed under a generalized Suslin
operation (the similar result for the Ramsey property was proved
by Matet). Further we compare two notions of forcing, the Mathias
forcing and a dual form of it, and will give some symmetries between
them. Finally we give some relationships between the dual-Mathias
forcing and the dual-Ramsey property.
\end{abstract}

\bigskip

\section{Notations and definitions}

Most of our set-theoretic notations and notations of
forcings are standard and can be found
in \cite{Jech} or \cite{Kunen}. An exception is that we will
write $A^B$ for the set of all functions from $B$ to $A$, instead
of\ \ ${}^B\hspace{-1.5mm}A$ because we never use ordinal arithmetic.
$A^{<\omega}$ is the set of all partial functions $f$
from $\omega$ to $A$ such that the cardinality of dom$(f)$ is finite.
\hfill\smallskip

First we will give the definitions of the sets we will consider as
the real numbers.\hfill\smallskip

Let $[x]^\kappa :=\{y\subseteq x: |y|=\kappa\}$ and
$[x]^{<\kappa} :=\{y\subseteq x: |y|<\kappa\}$,
where $|y|$ denotes the cardinality of $y$.
For $x\in [\omega]^\omega$, we will consider
$[x]^{<\omega}$ as the set of strictly increasing,
finite sequences in $x$
and $[x]^{\omega}$ as the set of strictly increasing,
infinite sequences in $x$. For $x\in\reals$ and $n\in\omega$ let $x(n)$
be such that $x(n)\in x$ and $|x(n)\cap x|=n$.\hfill\smallskip

We can consider $[\omega]^\omega$ also as
the set of infinite $0$-$1$-sequences (denoted by $2^{\omega}$)
or as the set of all infinite sequences in $\omega$ (denoted by
${\omega}^{\omega}$).

\subsubsection*{The Ellentuck topology}

We define a topology on $\reals$. Let
$X\in\reals$ and $s\in [\omega]^{<\omega}$ such that
$\mmax(s)<\mmin(X)$; then
$[s,X]^\omega:=\{Y\in\reals :Y\subs (s\cup X)\wedge s\subs Y\}$.
Now let the basic open sets on $\reals$ be the
sets $[s,X]^\omega$. These sets are called
the {\em Ellentuck neighborhoods.} The topology induced by the
Ellentuck neighborhoods is called the {\em Ellentuck topology\/}.

\subsubsection*{Relations on the set of partitions}

A {\em partition\/} $\pp$ (of $\omega$) is a subset of ${\cal P}
(\omega)$ such that the following holds:
\begin{tabbing}
iii)' \= \kill
i)\> if $b\in\pp$ then $b\neq\emptyset$\\
ii)\> if $b_1,b_2\in\pp$ and $b_1\neq b_2$ then $b_1\cap b_2 =\emptyset$\\
iii)\> $\bigcup \pp =\omega$.
\end{tabbing}

A partition means always a partition of $\omega$. If $\pp$ is a partition
and $b\in\pp$ then we call $b$ a block of $\pp$. If a partition has
infinitely many blocks (or equivalently if $\pp$ is infinite) we call
$\pp$ an infinite partition. The set of all infinite partitions is
denoted by $\parto$.\hfill\smallskip

If $\pp$ is a partition, $b\in\pp$ and $n,m\in\omega$ both
belong to $b$, then we write $\natural_{\pp}(n,m)$.
On the other hand with $\{\{n,m\}\in [\omega]^2:
\natural_{\pp}(n,m)\}$ we can reconstruct the
partition $\pp$.\hfill\smallskip

A {\em partial partition\/} $\pp'$ is a subset of ${\cal P}
(\omega)$ such that (i) and (ii) hold but instead of (iii) we have
\begin{tabbing}
iii)$\prime$ \= \kill
iii)${}'$\> $\bigcup \pp' =:\mdom(\pp')\subseteq\omega$.
\end{tabbing}

Note that a partition is always also a partial partition.
If $\mdom(\pp')\in\omega$ then $\pp'$ is a partition of
some $n\in\omega$.
The set of all partial partitions $\pp'$
where $\mdom(\pp')\in\omega$ is
denoted by $\NN$.
For $s\in\NN$, $s^*$ denotes the partial partition $s\cup\{\{
\mdom (s)\}\}$.\hfill\smallskip

Let $X_1,X_2$ be two partial partitions. We say that $X_1$ is {\em
coarser\/} than $X_2$, or that $X_2$ is {\em finer\/} than $X_1$, and
write $X_1\sqsubseteq X_2$ if for all blocks $b\in X_1$ the set
$b\cap\mbox{dom}(X_2)$ is the union of some sets
$b_i\cap\mbox{dom}(X_1)$, where each $b_i$ is a block of $X_2$. Let
$X_1\sqcap X_2$ denote the finest partial partition which is coarser
than $X_1$ and $X_2$ such that $\mbox{dom}(X_1\sqcap X_2)=
\mbox{dom}(X_1)\cup\mbox{dom}(X_2)$.\hfill\smallskip

If $f\in[\omega]^{<\omega}$ is a finite subset of $\omega$, then
$\{f\}$ is a partial partition with $\mbox{dom}(\{f\})=f$. For two partial
partitions $X_1$ and $X_2$ we write $X_1\sqsubseteq^* X_2$ if there is a
finite set $f\subseteq\mbox{dom}(X_1)$ such that
$X_1\sqcap\{f\}\sqsubseteq X_2$
and say that $X_1$ is coarser${}^*$ than $X_2$.
If $X_1\sqsubseteq^* X_2$, $X_2\sqsubseteq^* X_1$ and
$\mbox{dom}(X_1)=\mbox{dom}(X_2)$, then we write
$X_1\stackrel{*}{=} X_2$.\hfill\smallskip

Let $\pp_1,\pp_2$ be two partial partitions. If each block of $\pp_1$
can be written as the intersection of a block of $\pp_2$ with
$\mdom(\pp_1)$, then we write $\pp_1\seg\pp_2$.
Note that $\pp_1\seg\pp_2$ implies $\mdom(\pp_1)\subseteq\mdom(\pp_2)$.
\hfill\smallskip

If $\pp$ is a partial partition, then $\mMin (\pp)$\index{$\mMin ()$}
denotes the set $\{n\in\omega :\exists b\in\pp (n=\mmin(b))\}$,
where $\mmin(b):=\bigcap b$. If we order the blocks of $\pp$ by their
least element, then $\pp (n)$ denotes the $n$th block in this
ordering and $\pp (n)(k)$ denotes the $k$th element (in the natural
ordering) belonging to this block.

\subsubsection*{The dual Ellentuck topology}

We define a topology on the set of partitions as follows. Let
$X\in\parto$ and $s\in\NN$ such that $s\ceq X$. Then
${\open{s}{X}}:=\{Y\in\parto :s\seg Y\wedge Y\ceq X\}$ and $\partX
:={\open{\emptyset}{X}}$. Now let the basic open sets on $\parto$ be the
sets ${\open{s}{X}}$ (where $X$ and $s$ as above). These sets are called
the {\em dual Ellentuck neighborhoods.} The topology induced by the dual
Ellentuck neighborhoods is called the {\em dual Ellentuck topology\/}
(cf.~\cite{dual}).

\subsubsection*{Two notions of forcing}

The {\em Mathias forcing\/} $\M$ is defined as follows:
\[\begin{array}{c}
\langle s,S\rangle\in \M \Leftrightarrow s\in {[\omega ]}^{<\omega}
\ \wedge\ S\in {[\omega ]}^{\omega}\ \wedge\
\mbox{max($s$)$<$min($S$)},\\
\langle s,S\rangle\leq\langle t,T\rangle
\Leftrightarrow \ t\subs s\
\wedge\ S\subseteq T\ \wedge\ \forall n\in (s\setminus t)
(n\in T).
\end{array}\]

If $\langle s,S\rangle$ is an $\M$-condition, then we call $s$ the
{\em{stem}\/} of the condition. The Mathias forcing $\M$ has a lot of
combinatorial properties (which can be found in \cite{happy} and
\cite{delta1-2} or in \cite{Halb}). Note that we can consider an
$\M$-condition $\la s,S\ra$ as an Ellentuck neighborhood $[s,S]^\omega$
and $\la s,S\ra \leq\la t,T\ra$ if and only if $[s,S]^\omega\subs
[t,T]^\omega$.

The {\em dual-Mathias forcing\/} $\dM$ is
defined similarly to the Mathias forcing $\M$,
using the dual Ellentuck topology instead of the Ellentuck topology. So,
$$\langle s,X\rangle\in\dM\ \Leftrightarrow\ {\open{s}{X}}\ \mbox{is a dual
Ellentuck neighborhood}$$
and $$\langle s,X\rangle\leq\langle t,Y\rangle\ \Leftrightarrow\
{\open{s}{X}}\subseteq{\open{t}{Y}}.$$
If $\langle s,X\rangle$ is an $\dM$-condition, then we call $s$ again the
{\em{stem}\/} of the condition. Because the dual-Mathias forcing
is very close to the usual Mathias forcing, it also has some nice
properties similar to those of $\M$.

\subsubsection*{Two Ramsey properties}

The classical Ramsey property is a property of sets of infinite subsets
of $\omega$ (of sets of reals).
A set $A\subseteq\reals$ has
the {\em{Ramsey property}\/} (or is {\em{Ramsey}})
if $\exists X\in\reals([X]^\omega\subseteq A\vee
[X]^\omega\cap A=\emptyset).$
If there exists an $X$ such that
$[X]^\omega\cap A=\emptyset$ we call $A$ a
{\em{Ramsey null}} set. A set $A\subseteq\reals$ is
{\em{completely Ramsey}\/}
if for every Ellentuck neighborhood $[s,Y]^\omega$ there is an
$X\in [s,Y]^\omega$ such that $[s,X]^\omega\subs
A$ or $[s,X]^\omega\cap A=\emptyset$. If we are always in the latter
case, then we call $A$ {\em{completely Ramsey null.}}

The {\em{dual-Ramsey property}\/} deals with sets of
infinite partitions of
$\omega$.
A set $A\subseteq\parto$ has
the {\em{dual-Ramsey property}\/} (or is {\em{dual-Ramsey}})
if $\exists X\in\parto(\partX\subseteq A\vee
\partX\cap A=\emptyset).$
If there exists an $X$ such that
$\partX\cap A=\emptyset$ we call $A$ a
{\em{dual-Ramsey null}} set. A set $A\subseteq\parto$ is
{\em{completely dual-Ramsey}\/}
if for every dual Ellentuck neighborhood ${\open{s}{Y}}$ there is an
$X\in {\open{s}{Y}}$ such that
${\open{s}{X}}\subs A$
or ${\open{s}{X}}\cap A=\emptyset$. If we are always in the latter
case, then we call $A$ {\em{completely dual-Ramsey null.}}

Now we can start to give some symmetries between the two Ramsey properties
and between the two Mathias forcings.

\section{Basic facts}
\label{sec:basicfacts}

In this section we give the tools to consider sets of partitions as
sets of reals and to compare the two Ramsey properties. We will give
also some basic facts and well-known results concerning
the dual-Ramsey property and dual-Mathias forcing.
Further we give some symmetries
between Mathias forcing and the dual-Mathias forcing.

To compare the two Ramsey properties we first show that we can consider
each $A\subs\reals$ as a set of infinite partitions of $\omega$ and
vice versa. For this we define some arithmetical relations and functions.

Let $n,m\in\omega$ then
$\mbox{div}(n,m):=\mbox{max}(\{k\in\omega:k\cdot m\leq n\}$.
For $\{n,m\}\in[\omega]^2$ let $\flat\{n,m\}:=\frac{1}{2}
({\mbox{max}(\{n,m\})^2-\mbox{max}(\{n,m\})})+\mbox{min}(\{n,m\})$.
Consider $\flat\{n,m\}$ as undefined for $n=m$.

Let $x\in [\omega]^\omega$; then $\mbox{trans}(x)\subseteq\omega$ is such
that $n\not\in\mbox{trans}(x)$ {\em iff\/} there is a finite sequence $s$
of natural numbers of length $l+1$ such that
$$n=\flat\{s(0),s(l)\}\ \ \mbox{and}\ \ \forall k\in\{1,\ldots\!,l\}
(\flat\{s(k-1),s(k)\}\not\in x).$$
Note that $\mbox{trans}(x)\subseteq x$.
If $x\in [\omega]^\omega$, then we can consider $x$ as a partition with
$$\natural_x (n,m)\ \mbox{if and only if}\ n=m\ \mbox{or}
\ \flat\{n,m\}\not\in\mbox{trans}(x).$$
The {\em corresponding partition\/} of a real
$x\in [\omega]^\omega$ is denoted
by $\cp (x)$. Note that $\cp (x)\in\parto$ {\em{iff}\/} $\forall k\exists
n>k\forall m<n(\neg\natural_x (n,m))$ and further if
$y\subseteq x$, then $\cp (y)\ceq\cp (x)$.

A partition $\pp$ of $\omega$ we encode by a real $\pac (\pp)$ (the
{\em partition code\/} of $\pp$) as follows.
$$\pac (\pp):=\{k\in\omega :\exists n m(k=\flat\{n,m\}\wedge
\neg\natural_{\pp}(n,m)\}.$$
Note that if
$\pp_1\ceq\pp_2$ then $\pac (\pp_1)\subseteq\pac (\pp_2)$. With these
definitions we get the

\begin{fct}
\label{fact:finer}
The dual Ellentuck topology is
finer than the topology of the Baire space.
\end{fct}

\proof
Let $s\in \omega^{<\omega}$ and $U_s =\{f\in\omega^\omega :s\subset f\}$
be a basic open set in the Baire space $\omega^\omega$.
Because there is a bijection between $\omega^\omega$ and $\reals$, we
can write $U_s$ as a set $V_{s'}=\{r\in\reals :s'\subset r\wedge
\mmin(r\setminus s)>\mmax(s)\}$. Now
$\cp [V_{s'}]\cap\parto$ (where $\cp [V_{s'}]:=\{\cp (r): r\in
V_{s'}\}$) is open with respect to the dual Ellentuck topology.
Therefore the dual Ellentuck topology is
finer than the topology of the Baire space.\eor

\rmk A similar result is true for the Ellentuck topology
(cf.\,\cite{Ellentuck}).

\begin{fct}
\label{fct:RamseyBaire}
A set $C\subs\parto$ is completely dual-Ramsey if and only if $C$ has
the Baire property with respect to the dual Ellentuck topology and it
is completely dual-Ramsey null if and only if it is meager with
respect to the dual Ellentuck topology.
\end{fct}

\proof This is proved in \cite{dual}.
\eor

\rmk The analogous result is known for the Ramsey property with respect
to the Ellentuck topology (cf.\,\cite{Ellentuck}).

\subsubsection*{Some symmetries between the two Mathias forcings}

If $\mathfrak{g}$ is $\dM$-generic over $V$ and $\mathfrak{g}'\in
(\mathfrak{g})^{\omega}$, then also $\mathfrak{g}'$ is $\dM$-generic over $V$
(cf.~\cite{dual} Theorem 5.5). From this it follows immediately that
$\dM$ is proper and therefore
does not collapse $\aleph_1$. (For the definition
of properness consider e.g.~\cite{tools}.)

Further, for any $\dM$-condition $\langle s,X\rangle$ and any
sentence $\Phi$ of the forcing language $\dM$ there is an
$\dM$-condition $\langle s,Y\rangle\leq \langle s,X\rangle$ such
that $\langle s,Y\rangle\force_{\dM}\Phi$ or
$\langle s,Y\rangle\force_{\dM}\neg\Phi$ (cf.~\cite{dual}
Theorem 5.2). This property is called {\em{pure decision}}.

\rmk The similar results for Mathias forcing $\M$ can be found in
\cite{happy} (or in \cite{multiple}).\hfill\smallskip

We can write the dual-Mathias forcing as a two step iteration where
the first is the forcing notion $\fU$.

Let
$\fU$ be the partial order defined as follows:
$$p\in{\fU}\ \Leftrightarrow\
  p\in\parto,$$
$$p\leq q \ \Leftrightarrow\ p\ceq^* q.$$

We can also write the Mathias forcing as a two step iteration,
where the first step is the forcing notion $\bU$.
Let $\J :=[\omega]^{<\omega}$
be the ideal of finite sets and
let $\langle {\reals}/{\J},\leq\ \rangle =:
\bU$ be the partial order defined as follows.
$p\in{\bU}\ \Leftrightarrow\ p\in\reals$ and
$p\leq q \ \Leftrightarrow\ p\setminus q\in\J$
(this is $p\subs^* q$).

\begin{fct}\label{fct:o-closed}
The forcing notion $\fU$ is $\aleph_0$-closed and if
$\fD$ is $\fU$-generic over $V$, then $\mMin (\fD)$
is a Ramsey ultrafilter in $V[\fD]$.
\end{fct}

\proof
Let $X_1\geq X_2\geq\ldots$ be a decreasing sequence in $\fU$.
Choose a sequence
$f_i$ $(i\in\omega)$ of finite sets of natural numbers, such that
$X_{i+1}\puni \{f_i\}\ceq X_i$. Define $y_0:=X_0(0)$ and $y_n:=X_n(k)$
where $k:=3+\bigcup_{i<n}(\bigcup f_i)$. Now $Y:=\{y_i:i\in\omega\}
\cup (\omega\setminus\bigcup_{i\in\omega}y_i)$ is coarser${}^*$ than
each $X_i$ $(i\in\omega)$ and therefore $\fU$ is $\aleph_0$-closed.

Now we claim that the set $\{\mMin(\pp):\pp\in\fD\}$ is a Ramsey
ultrafilter in $V[\fD]$.
Remember that a forcing notion which is $\aleph_0$-closed adds
no new reals to $V$ (cf.\cite{Jech} Lemma 19.6).
Take a $\pi\in 2^{[\omega]^2}$ and a
$Y\in\parto$; then by the Ramsey Theorem (cf.\cite{Jech} Lemma 29.1) for
$\mMin(Y)\in\reals$ there exists an infinite $r\subseteq \mMin(Y)$ such
that $\pi$ is constant on $[r]^2$. Now let $X:=\{b:b\in Y\wedge b\cap
r\neq\emptyset\}\cup\bigcup\{b:b\in Y\wedge b\cap r=\emptyset\}$; then
$X\ceq Y$ and $\mMin(X)=r$. Thus
$H_{\pi}:=\{X\in\parto: \pi |_{[\mMin(X)]^2}\ \mbox{is constant}\}$
is dense in $\fU$ and hence $H_{\pi}\cap\fD\neq\emptyset$.
\eor

\rmk It is easy to see that the forcing notion
$\bU$ is $\aleph_0$-closed. Further we have that if
$D$ is $\bU$-generic over $V$, then $D$
is a Ramsey ultrafilter in $V[D]$.\hfill\smallskip

The forcing notion $\fU$ is stronger than the forcing notion $\bU$.

\begin{fct}
\label{fct:U-generic}
If $\fD$ is $\fU$-generic, then
the set $\{\mMin(\pp):\pp\in\fD\}$ is
$\bU$-generic.
\end{fct}

\proof Let $A\subseteq\reals$ be a maximal anti-chain in $\bU$,
i.e., $A$ is a maximal almost disjoint family. Then
the set $D_A :=\{\pp\in\fU:\exists a\in A(\mMin(\pp)\subseteq^*
a)\}$ is dense in $\fU$.
\eor

We define now the second step of the two step iteration.

Let $\fF\subs\parto$.
The partial order $\Pf$ is defined as follows.
\begin{center}
{$\langle s,X\rangle\in \Pf \Leftrightarrow s\in\NN
\ \wedge\ X\in\fF\ \wedge\ {\open{s}{X}}\
\mbox{is a dual Ellentuck neighborhood},$}\\
{$\langle s,X\rangle\leq\langle t,Y\rangle
\Leftrightarrow \ {\open{s}{X}}\subseteq{\open{t}{Y}}$.}
\end{center}

\rmk For ${\cal F}\subs\reals$ we can define the partial order
$\Pcf$ similarly.

\begin{fct}
\label{fct:fequi}
Let $\tilde{\fD}$ be the canonical $\fU$-name for the
$\fU$-generic object; then $${\fU}*\Padname \approx \dM.$$
\end{fct}

\proof
\[\begin{array}{ccl}
{\fU}*\Padname &\ =\ &\{\langle p,\langle \tilde{s},\tilde{X}\rangle
\rangle :p\in {\fU}\wedge p{\force}_{\fU}
\langle \tilde{s},\tilde{X}\rangle\in\Padname\}\\
 &=&\{\langle p,\langle \tilde{s},\tilde{X}\rangle
  \rangle : p\in\parto\wedge p{\force}_{\fU}
 (\tilde{X}\in\tilde{\fD}\wedge
 \tilde{s}\ceq\tilde{X})\}.
\end{array}\]

Now the embedding
\[\begin{array}{rccl}
h:\ \ &\dM&\ \longrightarrow\ &\ {\fU}*\Padname\\
      &\langle s,X\rangle& \longmapsto &\langle X,
      \langle\check{s},\check{X}
      \rangle\rangle
\end{array}\]
is a dense embedding (see \cite{tools} Definition 0.8):

\begin{enumerate}
\item It is easy to see that $h$ preserves the order relation $\leq$.
\item Let $\langle p,\langle \tilde{s},\tilde{X}\rangle\rangle\in{\fU}*
   \Padname$.
   Because ${\fU}$ is $\aleph_0$-closed, there is a condition
   $q\leq p$ and $s\in\NN,\ X\in\parto$
   such that $q{\force}_{\fU}\check{s}=\tilde{s}\wedge
   \check{X}=\tilde{X}$.
   Evidently, $\langle q,\langle\check{s},\check{X}\rangle\rangle\in
   {\fU}*\Padname$ is stronger than $\langle p,\langle
   \tilde{s},\tilde{X}\rangle\rangle$. Let $Z:=q\puni X$
   and let $Z'\ceq^* Z$ be such that $s\ceq Z'$. Now we have
   $h(\langle s,Z'\rangle)\leq \langle p,\langle
   \tilde{s},\tilde{X}\rangle\rangle$. \eor
\end{enumerate}

\rmk Let $\tilde{D}$ be the canonical $\bU$-name for the
$\bU$-generic object, then
${\bU}*{\operatorname{\mathbf{P}}}_{\tilde{D}}\approx \M.$
\hfill\smallskip

The dual-Mathias forcing is stronger than the Mathias forcing.

\begin{fct}
\label{fct:add-Mathias}
The dual-Mathias forcing adds Mathias reals.
\end{fct}

\proof Let $\DD$ be ${\fU}$-generic over $V$;
then by Fact~\ref{fct:U-generic},
$D:=\{\mMin(X):X\in\DD\}$ is ${\bU}$-generic over $V$. Now we
define $h:\Pad\rightarrow\Pd$ as follows.
\[\begin{array}{rccl}
h:\ \ &\Pad&\ \longrightarrow\ &\ \Pd\\
      &\langle s,X\rangle& \longmapsto &\langle \mMin (s),
      \mMin (X)\setminus\mMin (s)\rangle
\end{array}\]
For $h$ the following is true.\hfill\\
(i) If $q_1,q_2\in\Pad$, $q_1\leq q_2$, then $h(q_1)\leq h(q_2)$.\\
(ii) $\forall q\in\Pad\forall p'\leq h(q)\exists q'\in\Pad$ such that $q$
and $q'$ are compatible and $h(q')\leq p'$.\hfill\\
Therefore with \cite{multiple} Part~I, Lemma~2.7 we finally get
$V^{\M}\subseteq V^{\dM}$.
\eor

\section{On the dual-Ramsey property}
\label{sec:dualRam}

In this section we will show that the dual-Ramsey property is closed
under a generalized Suslin operation. As a corollary we will get
the already known result that analytic sets are completely dual-Ramsey.

Let $\fJ\subseteq {\cal P}(\parto )$ be the set of all completely
dual-Ramsey null sets. Further let
$\add (\fJ)$ be the smallest cardinal $\kappa$ such that there
exists a family ${\cal F}=\{J_\alpha\in\fJ :\alpha <\kappa\}$ with
$\bigcup{\cal F}\not\in\fJ$ and let
$\cov (\fJ)$ be the smallest cardinal $\kappa$ such that there
exists a family ${\cal F}=\{J_\alpha\in\fJ :\alpha <\kappa\}$ with
$\bigcup{\cal F}=\parto$. In \cite{Lori} it is shown that
$\cov (\fJ)=\add (\fJ)= \fH$, where $\fH$ is the dual-shattering
cardinal. Further it is shown that $\fH >\omega_1$ is relatively
consistent with ZFC.\\
Let $\Seq (\kappa ):=\kappa^{<\omega}$ and for $f\in\kappa^{\omega}$,
$n\in\omega$, let $\bar{f}(n)$ denote the finite sequence $\langle f(0),
f(1),\ldots\!,\linebreak[2] f(n-1)\rangle$. The generalized Suslin operation
${\cal A}_{\kappa}$ (for a cardinal $\kappa$) is defined as follows:
$${\cal A}_{\kappa}\{Q_s: s\in\Seq (\kappa )\}:=
\bigcup\limits_{f\in\kappa^{\omega}}\bigcap\limits_{n\in\omega}
Q_{\bar{f}(n)}\,.$$
In Theorem~\ref{thm:suslin} below
we will show that for each cardinal $\kappa <\fH$,
the completely dual-Ramsey
sets are closed under the operation ${\cal A}_{\kappa}$. But first
we give some other results.
\hfill\smallskip

A set $R\subseteq\parto$ is {\em dual Ellentuck meager\/}
if $R$ is meager with
respect to the dual Ellentuck topology. Remember that
a set is dual Ellentuck meager if and only if it is
completely dual-Ramsey null and a set is completely dual-Ramsey
if and only if it has the Baire property with respect to the dual
Ellentuck topology.\\
If ${\open{s}{X}}$ is a dual Ellentuck neighborhood then we
say that $R$ is dual Ellentuck meager in
${\open{s}{X}}$ if $R\cap {\open{s}{X}}$ is dual Ellentuck meager. By
\cite{dual} Theorem~4.1, $R$ is dual Ellentuck meager in ${\open{s}{X}}$
if for all ${\open{t}{Y}}\subseteq {\open{s}{X}}$ there exists a
partition $Z\in{\open{t}{Y}}$ such that ${\open{t}{Z}}\cap
R=\emptyset$.\hfill\smallskip

Let $R\subseteq\parto$ and $M:=\bigcup\{
{\open{s}{X}}: R$ is dual Ellentuck meager in
${\open{s}{X}}\}$. Further let $M(R):=M\cap R$. We first show that
\begin{lm}\label{lm:R}
If ${\open{s}{X}}$ is a dual Ellentuck neighborhood such that
${\open{s}{X}}\subseteq M$, then $R$ is dual Ellentuck meager in
${\open{s}{X}}$.
\end{lm}

\proof If ${\open{s}{X}}\subseteq M$, then ${\open{s}{X}}=
\bigcup\{{\open{t}{Y}}\subseteq {\open{s}{X}}: R$ is dual Ellentuck
meager in ${\open{t}{Y}}\}$. Let $N:=\bigcup\{{\open{u}{Z}}\subseteq
{\open{s}{X}}: R\cap{\open{u}{Z}}=\emptyset\}$. Because $N$ is
an open set, $N$ is completely dual-Ramsey. Therefore, for any
${\open{t}{Y}}\subseteq {\open{s}{X}}$ there exists a partition
$Y' \in{\open{t}{Y}}$ such that ${\open{t}{Y'}}\subseteq N$ or
${\open{t}{Y'}}\cap N=\emptyset$. If we are in the latter case,
then because ${\open{t}{Y'}}\subseteq {\open{s}{X}}$, we find
a ${\open{u}{Y''}}\subseteq {\open{t}{Y'}}$ such that $R$ is dual
Ellentuck meager in ${\open{u}{Y''}}$. Hence, there exists a
${\open{u}{Z}}\subseteq {\open{u}{Y''}}$ such that ${\open{u}{Z}}
\cap R=\emptyset$, which contradicts ${\open{t}{Y'}}\cap N=\emptyset$.
So we are always in the former case, which implies that $R$ is
dual Ellentuck meager in ${\open{s}{X}}$.
\eop

With this result, we can easily prove the following
\begin{lm}
\label{lm:MR}
The set $M(R)$ is dual Ellentuck meager.
\end{lm}

\proof Take a dual Ellentuck neighborhood ${\open{s}{X}}$ and let
$S:=\bigcup\{{\open{t}{Z}}\subseteq{\open{s}{X}}:R$ is dual Ellentuck
meager in ${\open{t}{Z}}\}$. Then $S$ as the union of open sets is open
and a subset of ${\open{s}{X}}$. Because ${\open{s}{X}}$ is also closed
(in the dual Ellentuck topology), the set $C:={\open{s}{X}}\setminus S$
is closed. By \cite{dual} Theorem~4.1 the sets $C$ and $S$ both are
completely dual-Ramsey. Therefore we find for every
${\open{s'}{X'}}\subseteq{\open{s}{X}}$ a partition $Y\in{\open{s'}{X'}}$
such that ${\open{s'}{Y}}\subseteq S$ or ${\open{s'}{Y}}\subseteq C$.
Now if ${\open{s'}{Y}}\subseteq S$, then by Lemma~\ref{lm:R},
$R$ is dual Ellentuck meager in ${\open{s'}{Y}}$ and
if ${\open{s'}{Y}}\subseteq C$,
then ${\open{s'}{Y}}\cap M(R)=\emptyset$.
To see this, assume there is an $H\in
M(R)\cap{\open{s'}{Y}}$. Because $H\in M(R)$ there exists a dual
Ellentuck neighborhood ${\open{t}{Z}}$ such that $H\in{\open{t}{Z}}$ and
$R$ is dual Ellentuck meager in ${\open{t}{Z}}$. Because
$H\in{\open{t}{Z}}$ and $H\in{\open{s'}{Y}}$ there is a dual Ellentuck
neighborhood ${\open{u}{U}}\subseteq{\open{t}{Z}}\cap{\open{s'}{Y}}$. But
with ${\open{u}{U}}\subseteq{\open{t}{Z}}$ it follows that $R$ is dual
Ellentuck meager in ${\open{u}{U}}$ and therefore
${\open{u}{U}}\subseteq S$, a contradiction to
${\open{u}{U}}\subseteq{\open{s'}{Y}}\subseteq C$.

Therefore, in both cases $M(R)$ is dual Ellentuck meager in
${\open{s'}{Y}}\subseteq{\open{s'}{X'}}$ and because ${\open{s}{X}}$
and ${\open{s'}{X'}}\subseteq{\open{s}{X}}$ were
arbitrary, the set $M(R)$ is dual Ellentuck meager in each dual Ellentuck
neighborhood. Hence, the set $M(R)$ is dual Ellentuck meager.
\eop

\begin{cor}
\label{cor:MR}
The set $R\cup (\parto\setminus M)$ has the dual Ellentuck Baire property.
\end{cor}
\proof Because $M$ is open, $\parto\setminus M$
is closed and $R\cup (\parto\setminus M)=(R\cap
M)\cup (\parto\setminus M)=M(R)\cup (\parto\setminus M)$ which is the
union of a meager set and a closed set and therefore has the dual
Ellentuck Baire property.
\eop

\begin{thm}
\label{thm:MR2}
If $R\subseteq\parto$, then we can construct a set $A\supseteq R$ which
has the dual Ellentuck Baire property and
whenever $Z\subseteq A\setminus R$ has the dual Ellentuck Baire property,
then $Z$ is dual Ellentuck meager.
\end{thm}

\proof Let $A:=R\cup (\parto\setminus M)$ where $M:=\bigcup\{
{\open{s}{X}}: R$ is dual Ellentuck meager in ${\open{s}{X}}\}$. By
Lemma~\ref{lm:MR} and Corollary~\ref{cor:MR} we know that $A$ has the
dual Ellentuck Baire property. Now let $Z\subseteq A\setminus R$ with the
dual Ellentuck Baire property. If $Z$ is not dual Ellentuck meager, then
there exists a dual Ellentuck neighborhood ${\open{u}{U}}$, such that
${\open{u}{U}}\setminus Z$ and therefore ${\open{u}{U}}\cap R$
are dual Ellentuck meager. Hence, $R$ is dual Ellentuck meager in
${\open{u}{U}}$ and therefore ${\open{u}{U}}\subseteq M$.
Since ${\open{u}{U}}\cap Z\neq\emptyset$
and $Z\cap M=\emptyset$, there is a $Y\in{\open{u}{U}}$ such that
$Y\not\in M$, a contradiction to $R$ is dual Ellentuck meager in
${\open{u}{U}}$.
\eop

Now we can prove the following
\begin{thm}\label{thm:suslin}
Let $\kappa <\fH$ be a cardinal number and for each
$s\in\Seq (\kappa )$ let $Q_s\subseteq\parto$. If all the sets
$Q_s$ are completely dual-Ramsey, then the set
$${\cal A}_{\kappa}\{Q_s:s\in\Seq (\kappa )\}$$
is completely dual-Ramsey too.
\end{thm}

\proof
Let $\{Q_s:s\in\Seq (\kappa )\}$ be a set of completely dual-Ramsey sets
and let $A:={\cal A}_{\kappa}\{Q_s:s\in\Seq (\kappa )\}$.
For two sequences $s,f\in\kappa^{\leq\omega}$ we write $s\subseteq f$
if $s$ is an initial segment of $f$. If $s\in\kappa^{<\omega}$ is a
finite sequence, then $|s|$ denotes the length of $s$.
Without loss of generality
we may assume that $Q_s\supseteq Q_t$ whenever $s\subseteq t$.\\
For $s\in\Seq (\kappa )$ let
$$A_s:=\bigcup\begin{Sb}
f\in\kappa^{\omega} \\ s\subseteq f
\end{Sb}
\bigcap\begin{Sb}
n\in\omega \\ n\geq |s|
\end{Sb}
Q_{\bar{f}(n)}.$$
For $s\in\Seq (\kappa )$ we have $A_s\subseteq Q_s$,
$A_s=\bigcup_{\alpha <\kappa}
A_{s{\join{\alpha}}}$ and $A=A_{\emptyset}$.
By Theorem~\ref{thm:MR2}, for each $s\in\Seq (\kappa )$ we find
a $B_s\supseteq A_s$ which is completely dual-Ramsey and if
$Z\subseteq B_s\setminus A_s$ has the dual-Ramsey property, then
$Z$ is dual-Ramsey null. Because $Q_s\supseteq A_s$ is completely
dual-Ramsey, we may assume that $B_s\subseteq Q_s$ and therefore
$$A={\cal A}_{\kappa}\{B_s:s\in\Seq (\kappa )\}.$$
Let $B:=B_{\emptyset}$. Note that $A=\bigcup_{\alpha<\kappa}
A_{\langle\alpha\rangle}\subseteq\bigcup_{\alpha<\kappa}
B_{\langle\alpha\rangle}$ and therefore $B\subseteq
\bigcup_{\alpha<\kappa} B_{\langle\alpha\rangle}$. Now we show that
$$B\setminus A\subseteq \bigcup\limits_{\alpha<\kappa}
B_{\langle\alpha\rangle}\subseteq \bigcup
\limits_{f\in\kappa^{\omega}}\bigcap\limits_{n\in\omega}
B_{\bar{f}(n)}\subseteq\bigcup\limits_{s\in\Seq(\kappa)}(B_s
\setminus\bigcup\limits_{\alpha<\kappa}B_{s{\join{\alpha}}})\,.$$
Assume $x\not\in\bigcup_{s}(B_s
\setminus\bigcup_{\alpha<\kappa}B_{s{\join{\alpha}}})$.
If we have for all $\alpha<\kappa$, that
$x\not\in B_{\langle\alpha\rangle}$, then $x\not\in B$.
And if there exists an $\alpha_0<\kappa$ such that $x\in B_{\langle
\alpha_0\rangle}$, because $x\not\in\bigcup_{s}
(B_s\setminus\bigcup_{\alpha<\kappa}B_{s{\join{\alpha}}})$ we find an
$\alpha_1$ such that $x\in B_{\langle\alpha_0,\alpha_1\rangle}$
and finally we find an $f\in\kappa^{\omega}$ such that for all
$n\leq\omega$: $x\in B_{\bar{f}(n)}$.
But this implies that $x\in A$.
Now because
$B_s\setminus\bigcup_{\alpha<\kappa}B_{s{\join{\alpha}}}
\subseteq B_s\setminus\bigcup_{\alpha<\kappa}A_{s{\join{\alpha}}}
=B_s\setminus A_s$ and because
$\bigcup_{\alpha<\kappa}B_{s{\join{\alpha}}}$ is the union
of less than $\fH$ completely dual-Ramsey sets,
$B_s\setminus\bigcup_{\alpha<\kappa}B_{s{\join{\alpha}}}$ is
completely dual-Ramsey and as a subset of $B_s\setminus A_s$, it
is completely dual-Ramsey null. Therefore, $B\setminus A$ as a subset
of the union of less than $\fH$ completely dual-Ramsey null sets
is completely dual-Ramsey null and because $B$ is completely
dual-Ramsey, $A$ is completely dual-Ramsey too.
\eop

\rmk A similar result holds also for the Ramsey property and is
proved by Matet in \cite{Matet}.\hfill\smallskip

As a corollary we get a result which was first proved by Carlson
and Simpson (cf.\,\cite{dual}).

\begin{cor}
\label{cor:analytic}
Every analytic set is completely dual-Ramsey.
\end{cor}

\proof This follows from Theorem~\ref{thm:suslin} and because
each analytic set $A\subseteq\reals$ can be written as
$$A={\cal A}\{Q_s: s\in\Seq (\omega )\}$$
where each $Q_s\subseteq\reals$ is a closed set in the Baire space.
\eop

\rmk For a similar result cf.\,\cite{Ellentuck} or \cite{Silver}.

\section{Game-families and the forcing notion P${}_{\fF}$}
\label{sec:gamefamilies}

First we define a game and game-families. Then we show that, for
game-families $\fF$, the forcing notion $\Pf$ has pure
decision and if $X$ is $\Pf$-generic and $Y\in (X)^{\omega}$, then
$Y$ is $\Pf$-generic too.

We call a family $\fF\subseteq\parto$ {\em non-principal}, if for
all $X\in\fF$ there is a $Y\in\fF$ such that $Y\ceq X$ and
$\neg (Y\ass X)$. A family $\fF$
is {\em closed under refinement}, if $X\ceq Y$ and $X\in\fF$ implies
that $Y\in\fF$. Further it is {\em closed under finite changes\/} if
for all $s\in\NN$ and $X\in\fF$, $X\kap s\in\fF$.

In the sequel $\fF$ is always a non-principal family which is closed
under refinement and finite changes.

If $s\in\NN$ and $s\ceq X\in\fF$, then we call the dual Ellentuck
neighborhood ${\open{s}{X}}$ an $\fF$-dual Ellentuck neighborhood
and write ${\opf{s}{X}}$ to emphasize that $X\in\fF$.
A set $\cO\subs\parto$ is called $\fF$-open if $\cO$ can be written
as the union of $\fF$-dual Ellentuck neighborhoods.
%A set $A\subs
%\parto$ is called $\fF${\em -dual-Ramsey}, if for every $\fF$-dual
%Ellentuck neighborhood ${\opf{s}{X}}$ there exists a $Y\in
%{\opf{s}{X}}\cap\fF$ such that ${\opf{s}{Y}}\subs A$ or ${\opf{s}{X}}
%\cap A\cap \fF =\emptyset$.
For $s\in\NN$ remember that $s^* =s\cup\{\{\mdom (s)\}\}$.
%If $s\ceq X\in\parto$, then
%$r_{s,X}\in\NN$ is such that $s\seg r_{s,X}$, $|s|=|r_{s,X}|$ and
%$r_{s,X}^*\seg X$. ($r_{s,X}$ is the root of ${\open{s}{X}}$.)
\hfill\smallskip

Fix a family $\fF\subs\parto$ (which is non-principal and closed under
refinement and finite changes). Let $X\in\fF$ and $s\in\NN$ be such
that $s\ceq X$.  We associate with ${\opf{s}{X}}$ the following game.
(This type of game was suggested first by Kastanas in \cite{Kastanas}.)

$$\begin{array}{lcccccccc}
{\operatorname{I}} &\ \ \ & \la X_0\ra &  &\la X_1\ra & &\la X_2\ra & &  \\
   &      &            &   &           &   & &          & \ldots \\
{\operatorname{II}}& & &\la t_0,Y_0\ra & &\la t_1,Y_1\ra & &\la t_2,Y_2\ra &
\end{array}$$

All the $X_i$ of player I and the $Y_i$ of player II must be elements of
the family $\fF$.
Player I plays $\la X_0\ra$ such that $X_0\in {\opf{s}{X}}$,
then II plays $\la t_0,Y_0\ra$ such that $Y_0\in {\opf{s}{X_0}}$,
$s\seg t_0^*\seg Y_0$ and $|t_0|=|s|$. For $n\geq 1$, the $n$th
move of player I is $\la X_n \ra$ such that $X_n\in
{\opf{t_{n-1}^*}{Y_{n-1}}}$ and then player II plays $\la t_n,Y_n\ra$
such that $Y_n\in {\opf{t_{n-1}^*}{X_{n}}}$, $t_{n-1}^*\seg t_n^*\seg Y_n$
and $|t_n|=|t_{n-1}|+1$. Player I wins iff the only $Y$ with $t_n\seg Y$
(for all $n$) is in $\fF$.
We denote this game by $\GF$ {\em {starting with}} $\la s,X\ra$.
\hfill\smallskip

A non-principal family $\fF$ which is closed under refinement and
finite changes is a {\em game-family} if player II has no winning
strategy in the game $\GF$.\hfill\smallskip

A family $\fF\subs\parto$ is called a {\em filter} if for any
$X,Y\in\fF$, also $X\kap Y\in\fF$. A filter which is also a game-family
is called a {\em game-filter}. Note that $\parto$ is game-family but
not a game-filter.
(But it is consistent with ZFC that game-filters exist, as
Theorem~\ref{thm:gamefilter} will show).

Let $\cO\subs\parto$ be an $\fF$-open set.
Call ${\opf{s}{X}}$ {\em good\/} (with respect to $\cO$),
if for some $Y\in {\opf{s}{X}}\cfF$,
${\opf{s}{Y}}\subs\cO$; otherwise call it {\em bad}. Note that
if ${\opf{s}{X}}$ is bad and $Y\in {\opf{s}{X}}\cfF$, then
${\opf{s}{Y}}$ is bad, too.
We call ${\opf{s}{X}}$ {\em ugly}
if ${\opf{t^*}{X}}$ is bad for all $s\seg t^*\ceq X$ with $|t|=|s|$.
Note that if ${\opf{s}{X}}$ is ugly, then ${\opf{s}{X}}$ is bad, too.
\hfill\smallskip

To prove the following two lemmas, we will follow in fact the proof
of Lemma~19.15 in \cite{Kechris}.

\begin{lm}\label{lm:ugly}
Let $\fF$ be a game-family and $\cO\subs\parto$ an $\fF$-open set.
If ${\opf{s}{X}}$ is bad (with respect to $\cO$), then there exists
a $Z\in {\opf{s}{X}}$ such that ${\opf{s}{Z}}$ is ugly.
\end{lm}

\proof
We begin by describing a strategy for player II in
the game $\GF$ starting
with $\la s,X\ra$. Let $\la X_n\ra$ be the $n$th move of player I
and $t_n$ be such that $s\seg t_{n}$, $|t_{n}|=|s|+n$
and $t_{n}^*\seg X_n$.
Let $\{t_{n}^i:i\leq m\}$ be an
enumeration of all $t$ such that $s\seg t\ceq
t_{n}$, $|t|=|s|$ and $\mdom (t)=\mdom (t_{n})$. Further
let $Y^{-1}:=X_{n}$. Now choose for each $i\leq m$ a partition $Y^i\in\fF$
such that $Y^i\ceq Y^{i-1}$, $t_{n}^*\seg Y^i$
and ${\opf{(t_{n}^i)^*}{Y^i}}$ is bad or
${\opf{(t_{n}^i)^*}{Y^i}}\subs\cO$.
Finally, let $Y_{n}:=Y^m$ and let player II play $\la t_{n},Y_{n}\ra$.
\hfill\smallskip

Because player II has no winning strategy, player I can play so
that the only $Y$ with $t_n\seg Y$ (for all $n$) belongs to $\fF$.
Let $S_Y:=\{t^*\ceq Y: s\seg t\wedge |t|=|s|\}$; then
(because of the strategy of player II), for all
$t\in S_Y$ we have either ${\opf{t^*}{Y}}$ is bad or ${\opf{t^*}{Y}}
\subs\cO$. Now let $C_0:=\{t\in S_Y: {\opf{t}{Y}}$ is bad$\}$
and $C_1:=\{t\in S_Y:{\opf{t^*}{Y}}\subs\cO\}=S_Y\setminus C_0$.
By a result of \cite{HalMat}, there exists a partition
$Z\in {\opf{s}{Y}}\cfF$, such that $S_Z\subs C_0$ or $S_Z\subs C_1$.
If we are in the latter case,
we have ${\opf{s}{Z}}\subs\cO$, which contradicts that
${\opf{s}{X}}$ is bad. So we must have $S_Z\subs C_0$, which implies
that ${\opf{s}{Z}}$ is ugly
and completes the proof of the Lemma.
%\ref{lm:ugly}
\eop

\begin{lm}\label{lm:2ndstep}
If $\fF$ is a game-family and $\cO\subs\parto$ is an $\fF$-open set, then
for every $\fF$-dual Ellentuck neighborhood ${\opf{s}{X}}$ there exists
a $Y\in {\opf{s}{X}}\cap\fF$ such that ${\opf{s}{Y}}\subs\cO$ or
${\opf{s}{Y}}\cap\cO\cap\fF =\emptyset$.
\end{lm}

\proof
If ${\opf{s}{X}}$ is good, then we are done. Otherwise we consider the
game $\GF$ starting
with $\la s,X\ra$. Let $\la X_0\ra$ be the first move of player I.
Because ${\opf{s}{X_0}}$ is bad, by Lemma~\ref{lm:ugly} we can choose $Y'
\in {\opf{s}{X_0}}\cfF$ such that ${\opf{s}{Y'}}$ is ugly.
Let $t_0$ be such that
$s\seg t_0^*\seg Y'$ and $|t_0|=|s|$. Now we choose $Y_0
\in {\opf{t_0^*}{Y'}}\cfF$ such that ${\opf{t_0^*}{Y_0}}$ is ugly,
which is is possible because ${\opf{t_0}{Y'}}$ is ugly and therefore
${\opf{t_0^*}{Y'}}$ is bad.
Note that for all $t$ with
$s\seg t\ceq t_{0}$ and $\mdom (t)=\mdom (t_{0})$ we have
${\opf{t^*}{Y_0}}$ is ugly. Now player II plays $\la t_0,Y_0\ra.$ \\
Let $\la X_{n+1}\ra$ be the $(n+1)$th move of player I. By the strategy
of player II we have ${\opf{t^*}{X_{n+1}}}$ is ugly
for all $t$ with $s\seg t\ceq t_{n}$ and $\mdom (t)=\mdom (t_{n})$.
Let $t_{n+1}$ be such that $|t_{n+1}|=|t_n|+1=|s|+n$ and
$t_n^*\seg t_{n+1}^*\seg X_{n+1}$.
Let $\{t_{n+1}^i:i\leq m\}$ be an
enumeration of all $t$ such that $s\seg t\ceq
t_{n+1}$ and $\mdom (t)=\mdom (t_{n+1})$. Further
let $Y^{-1}:=X_{n+1}$. Now choose for each $i\leq m$ a partition $Y^i\in\fF$
such that $Y^i\ceq Y^{i-1}$, $t_{n+1}^*\seg Y^i$
and ${\opf{(t_{n+1}^i)^*}{Y^i}}$ is ugly. (This is possible because we know
that ${\opf{t^*}{X_{k}}}$ is ugly for all $k\leq n$ and $t$ with
$s\seg t\ceq t_{k}$ and $\mdom (t)=\mdom (t_{k})$, which implies that
${\opf{(t_{n+1}^i)^*}{X_{n+1}}}$ is bad.)
Finally, let $Y_{n+1}:=Y^m$ and let player II play $\la t_{n+1},Y_{n+1}\ra$.
\hfill\smallskip

Because player II has no winning strategy, player I can play so
that the only $Y$ with $t_n\seg Y$ (for all $n$) belongs to $\fF$.
We claim that ${\opf{s}{Y}}\cap\cO\cap\fF =\emptyset.$ Let $Z\in
{\opf{s}{Y}}\cap\cO\cap\fF$. Because $\cO$ is $\fF$-open we find
a $t\seg Z$ such that ${\opf{t^*}{Z}}\subs\cO$. Because $t^*\ceq Y$
we know by the strategy of player II that ${\opf{t^*}{Y}}$ is bad.
Hence, there is no $Z\in {\opf{t^*}{Y}}$ such that
${\opf{t^*}{Z}}\subs\cO$. This completes the proof.
\eop

Now we give two properties of the forcing notion $\Pf$, (where
$\Pf$ is defined as in section~\ref{sec:basicfacts} and $\fF$ is
a game-family). Note that for $\fF=\parto$ (which is obviously
a game-family) the forcing notion $\Pf$ is the same as dual-Mathias
forcing.
The first property of the forcing notion $\Pf$ we give
is called {\em pure decision.}

\begin{thm}\label{thm:pd}
Let $\fF$ be a game-family and let $\Phi$ be a sentence
of the forcing language $\Pf$. For any $\Pf$-condition
${\opf{s}{X}}$ there exists a $\Pf$-condition
${\opf{s}{Y}}\leq {\opf{s}{X}}$ such that
${\opf{s}{Y}}\force_{\Pf} \Phi$ or
${\opf{s}{Y}}\force_{\Pf} \neg\Phi$.
\end{thm}

\proof With respect to $\Phi$ we define
$\cO_1:=\{Y:{\opf{t}{Y}}{\force}_{\Pf} \Phi\ \mbox{for some $t\seg
Y\in\fF$}\}$ and $\cO_2:=\{Y:{\opf{t}{Y}}{\force}_{\Pf}
\neg\Phi\ \mbox{for some $t\seg Y\in\fF$}\}$. Clearly $\cO_1$ and
$\cO_2$ are both $\fF$-open and $\cO_1\cup\cO_2$ is even dense (with
respect to the partial order $\Pf$).
Because $\fF$ is a game-family, by Lemma~\ref{lm:2ndstep} we know
that for any ${\opf{s}{X}}\in\Pad$ there exists
$Y\in{\opf{s}{X}}\cfF$ such that either ${\opf{s}{Y}}
\subs\cO_1$ or ${\opf{s}{Y}}\cap\cO_1\cap\fF =\emptyset$.
In the former case
we have ${\opf{s}{Y}}{\force}_{\Pf}\Phi$ and we are done.
In the latter case we find $Y'\in {\opf{s}{Y}}\cfF$ such that
${\opf{s}{Y'}}\subs\cO_2$. (Otherwise we would have ${\opf{s}{Y'}}\cap
(\cO_2\cup\cO_1)\cap\fF =\emptyset$,
which is impossible by the density of $\cO_1\cup\cO_2$.)
Hence, ${\opf{s}{Y'}}{\force}_{\Pf}\neg\Phi$.
\eop

Let $\fF$ be a game-family,
$\fG$ be $\Pf$-generic and define $X_\fG:=\bigcap\fG$. Now $X_\fG$
is an infinite partition and $\fG=\{{\opf{s}{Z}}: s\seg X_\fG\ceq Z\}$.
Therefore we can consider the partition
$X_\fG\in\parto$ as a $\Pf$-generic object.
Further we have $\fG\subs\Pf$ is $\Pf$-generic if and only if
$X_\fG \in\bigcup D$ for all $D\subs\Pf$ which are dense in $\Pf$.
Note that if $D$ is dense in $\Pf$, then $\bigcup D$ is $\fF$-open.
\hfill\smallskip

The next theorem shows in fact that if $\fF$ is a game-family, then $\Pf$ is
proper.

\begin{thm}\label{thm:proper}
Let $\fF\subs\parto$ be a game-family. If $X_0\in\parto$ is $\Pf$-generic
over $V$ and $Y_0\in (X_0)^{\omega}\cap V[X_0]$, then $Y_0$ is also
$\Pf$-generic over $V$.
\end{thm}

\proof
Take an arbitrary dense set
$D\subs\Pf$, i.e. for all ${\opf{s}{X}}$ there
exists a ${\opf{t}{Y}}\subs {\opf{s}{X}}$ such that ${\opf{t}{Y}}\in D$.
Let $D'$ be the set of all ${\opf{s}{Z}}$ such that ${\opf{t}{Z}}\subs
\bigcup D$ for all $t\ceq s$ with $\mdom (t)=\mdom (s)$.\hfill\smallskip

First we show that $D'$ is dense in $\Pf$.
For this take an arbitrary ${\opf{s}{W}}$
and let $\{t_i: 0\leq i\leq m\}$ be an enumeration of all $t\in\NN$
such that $t\ceq s$ and $\mdom (t)=\mdom (s)$. Because $D$ is dense
in $\Pf$ and $\bigcup D$ is $\fF$-open, we
find for every $t_i$ a $W'\in\fF$ such that $t_i\ceq W'$ and
${\opf{t_i}{W'}}\subs\bigcup D$. Moreover, if we define $W_{-1}:=W$, we
can choose for every $i\leq m$ a partition $W_i\in\fF$
such that $W_i\ceq W_{i-1}$,
$s\seg W_i$ and ${\opf{t_i}{W_i}}\subs\bigcup D$. Now
${\opf{s}{W_m}}\in D'$ and because ${\opf{s}{W_m}}\subs{\opf{s}{W}}$,
$D'$ is dense in $\Pf$.\hfill\smallskip

Since $D'$ is dense and $X_0\in\parto$ is $\Pf$-generic, there exists
a ${\opf{s}{Z}}\in D'$ such that $s\seg X_0\ceq Z$. Because $Y_0\in
(X_0)^{\omega}$ we have $t\seg Y_0\ceq Z$
for some $t\ceq s$ and because ${\opf{t}{Z}}\subs\bigcup D$, we get
$Y_0\in\bigcup D$. Hence, $Y_0\in\bigcup D$ for every dense $D\subs\Pf$,
which completes the proof.
\eop

\rmk A similar result is proved in \cite{happy} and \cite{Matet3}.

\section{On the dual-Mathias forcing and game-filters}
\label{sec:dMathias}

Now we show that it is consistent with ZFC that game-filters exist.
(Remember that a game-filter $\fF$ is a game-family which is also a
filter and a game-family is a non-principal family which is closed under
refinement and finite changes such that player II has no winning
strategy in the game $\GF$.) Further we show that the dual-Mathias
forcing $\dM$ is flexible and with this result we can prove that if
$V$ is $\Sigma^1_4$-$\dM$-absolute, then $\oV$ is inaccessible
in L.

In the sequel let $\fU$ be the forcing notion we defined
in section~\ref{sec:basicfacts}.

\begin{thm}\label{thm:gamefilter}
If $\fD$ is $\fU$-generic over $V$, then $\fD$ is a game-filter
in $V[\fD ]$ with respect to the game $\cG (\fD )$.
\end{thm}

\proof Because $\fD$ is $\fU$-generic over $V$, we know that $\fD
\subs\parto$ is a non-principal family in $V[\fD ]$ which is closed
under refinement and finite changes, and for $X,Y\in\fD$ we also
have $X\kap Y\in\fD$. It remains to show that player II has no winning
strategy in the game $\cG (\fD )$.\hfill\smallskip

Let $\tsig$ be a $\fU$-name for a strategy for player II
in the game $\cG (\nameDD )$,
where $\nameDD$ is the canonical $\fU$-name for the $\fU$-generic object.
Let us assume that player II will follow this strategy. We may assume
that $$\one\force_{\fU}``\tsig\ \mbox{is a strategy for II in the game}
\ \cG (\nameDD )".$$
If $$Z\force_{\fU}\tsig (\la\tX_0\ra,\la\tit_0,\tY_0\ra,\ldots,\la\tX_n\ra)
=\la\tit_n,\tY_n\ra,$$ then for $n\geq 1$ we get
$$Z\force_{\fU} (|\tit_n|=|\tit_{n-1}|+1\wedge\tit_{n-1}^*\seg\tit_{n}^*\seg
\tY_n\ceq\tX_n\wedge\tY_n\in\nameDD)$$
and for $n=0$ we have
$$Z\force_{\fU} (|\tit_0|=|\ts|\wedge\ts\seg\tit_{0}^*\seg
\tY_0\ceq\tX_0\ceq\tX\wedge\tY_0\in\nameDD)$$
where $\la\ts,\tX\ra$
is the starting point of $\cG (\nameDD )$.\hfill\smallskip

Now let $\la\ts,\tX\ra$ (the starting point of the game $\cG (\nameDD )$)
be such that ${\open{\ts}{\tX}}$ is a $\fU$-name for a dual Ellentuck
neighborhood and let $Z_0\in\parto\cap V$ be a $\fU$-condition in $V$
such that $Z_0\force_{\fU} \tX\in\nameDD$. Therefore,
$Z_0\force_{\fU}``{\open{\ts}{\tX}}$ is a $\nameDD$-dual Ellentuck
neighborhood". By Fact~\ref{fct:o-closed} we know that the forcing notion
$\fU$ adds no new reals (and therefore no new partitions) to $V$. So, we
find a $Z_0^{\prime}\ceq^* Z_0$ and a dual Ellentuck neighborhood
${\open{s}{X}}$ in $V$ such that
$$Z_0^{\prime}\force_{\fU}\la\ts,\tX\ra =\la\chs,\chX\ra$$
where $\chs$ and $\chX$ are the canonical $\fU$-names for $s$ and $X$.
Because $Z_0^{\prime}\force_{\fU}\chX\in\nameDD$, we must have
$Z_0^{\prime}\leq X$, which is the same as $Z_0^{\prime}\ceq^* X$.
Finally put $X_0\in\parto$ such that $X_0\ass Z_0^{\prime}$ and
$X_0\in {\open{s}{X}}$. Player I plays now $\la\chX_0\ra$.
Since player II follows the strategy $\tsig$, player II plays now
$\tsig (\la\chX_0\ra )=:\la\tit_0,\tY_0\ra$. Again by Fact~\ref{fct:o-closed}
there exists a
$Z_1\ceq^* X_0$ and a dual Ellentuck neighborhood ${\open{t_0}{Y_0}}$ in
$V$ such that
$$Z_1\force_{\fU}\la\tit_0,\tY_0\ra =\la\cht_0,\chY_0\ra.$$
And again by $Z_1\force_{\fU}\chY_0\in\nameDD$ we find $X_1\ass Z_1$
such that $t_0^*\seg X_1\ceq Y_0$. Player I plays now $\la\chX_1\ra$.
\hfill\smallskip

In general, if $\tsig (\la\tX_0\ra,\la\tit_0,\tY_0\ra,\ldots,\la\tX_n\ra)
=\la\tit_n,\tY_n\ra$, then player I can play $\chX_{n+1}$ such that
$X_n\force_{\fU}\la\tit_n,\tY_n\ra = \la\cht_n,\chY_n\ra$ and
$t_n^*\seg X_{n+1}\ceq Y_n$. For $n\geq m$ we also have $X_n\ceq X_m$.
Let $Y\in\parto$ be the such that $t_n\seg Y$ (for all $n$), then
$$Y\force_{\fU} ``\mbox{the only $\tY$ such that $\tit_n\seg\tY$ (for all
$n$) is in $\nameDD$}".$$
Hence, the strategy $\tsig$ is not a winning strategy for player II and
because $\tsig$ was an arbitrary strategy, player II has no winning
strategy at all.
\eop

\rmk A similar result is in fact proved in \cite{happy}
(cf. also \cite{Matet}).

As a corollary we get that the forcing notion
$\Pad$ (where ${\fD}$ is $\fU$-generic over $V$)
has pure decision in $V[\fD]$.

\begin{cor}\label{cor:pd}
Let $\fD$ be $\fU$-generic over $V$. Then the forcing notion
$\Pad$ has pure decision in $V[\fD ]$.
\end{cor}

\proof This follows from Theorem~\ref{thm:pd} and
Theorem~\ref{thm:gamefilter}.
\eop

Corollary~\ref{cor:pd} follows also from the facts that
the dual-Mathias forcing has pure decision (cf.\,\cite{dual})
and that it can be written as a two step iteration
as in section~\ref{sec:basicfacts}.\hfill\smallskip

\rmk If $D$ is $\bU$-generic over $V$, then $\Pd$ has pure decision in
$V[D]$ (cf.\,\cite{happy}).

\subsubsection*{Some more properties of $\dM$}

Let $\Pp$ be a notion of forcing in the model $V$. We say that $V$ is
$\Sigma^1_n$-$\Pp$-absolute if for every $\Sigma^1_n$-sentences $\Phi$
with parameters in $V$ the following is true.
$$V\models\Phi\ \mbox{if and only if}\ V[G]\models\Phi,$$
where $G$ is any $\Pp$-generic object over $V$.

Now we will show
that if $V$ is $\Sigma^1_4$-$\dM$-absolute, then
$\oV$ is inaccessible in $\LL$. For this we first will translate
the dual-Mathias forcing in a tree forcing notion.

If $s$ is a partial partition of some natural number $n\in\omega$, then
we can consider $s$ as a subset of ${\cal P}(n)$ or equivalently, as a
finite set of finite sets of natural numbers. Let $t$ be a finite set of
natural numbers, then $\sharp t$ is such that for all $k\in\omega:$
$\mbox{div}(\sharp t,2^k)$ is odd $\Leftrightarrow k\in s$. (Remember
that $\mbox{div}(n,m):=\mbox{max}(\{k\in\omega:k\cdot m\leq n\}$.)
Now $\sharp s$ is such that for all $k\in\omega:$
$\mbox{div}(\sharp s,2^k)$ is odd $\Leftrightarrow k=\sharp t$ for
some $t\in s$. (In fact $\sharp s$ is defined for any finite set of
finite sets of natural numbers.) If $s\in\NN$,
then $|s|$ denotes the cardinality of $s$, which is the number of blocks of
$s$.\hfill\smallskip

For $s\in\NN$ with $|s|=k$
let $\bar{s}$ be the finite sequence $\langle n_1,\ldots,n_k\rangle$
where $n_i:=\sharp s_i$ and $s_i\in\NN$ is such that $|s_i|=i$ and
$s_i^*\seg s^*$.

Now let $p={\open{s}{X}}$ be an
$\dM$-condition. Without loss of generality
we may assume that $s^*\ceq X$.
The tree ${\mathfrak{t}}_p\subseteq\omega^{<\omega}$ is defined as follows.
$$\sigma\in {\mathfrak{t}}_p\ \Leftrightarrow\ \exists t\in\NN ((t^*\seg s^*
\vee s\seg t)\wedge t^*\ceq X \wedge \sigma = \bar{t}.$$
\begin{fct}
\label{fct:easy}
Let $p,q$ be two $\dM$-conditions. Then ${\mathfrak{t}}_p$ is a
subtree of ${\mathfrak{t}}_q$ if and only if $p\leq q$.\eor
\end{fct}

Finally let $T_{\dM}:= \{ {\mathfrak{t}}_p :p\in
\dM \}$; then $T_{\dM}$ is
a set of trees. We stipulate that ${\mathfrak{t}}_p\leq {\mathfrak{t}}_q$
if ${\mathfrak{t}}_p$ is a subtree of ${\mathfrak{t}}_q$. Then (by
Fact~\ref{fct:easy}) forcing with ${\fT}_{\dM}:=\langle
T_{\dM},\leq\rangle$ is the same as forcing with $\dM$.

Now we will give the definition of a flexible forcing notion $\Pp$. But
first we have to give some other definitions.

A set $T\subseteq {\omega}^{<\omega}$ is
called a {\em{Laver-tree}\/} if
$$T\ \mbox{is a tree and}\ \exists\tau\in T\forall\sigma\in T
(\sigma\subseteq\tau \vee (\tau\subseteq\sigma\ \wedge\
|\{n:\sigma^{\frown}n\in T\}| =\omega)).$$
(We call $\tau$ the stem of $T$. For $\sigma\in T$ we let succ${}_T
(\sigma ):=\{n:\sigma^{\frown}n\in T\}$, (the successors of $\sigma$
in $T$) and $T_\rho :=\{\sigma\in T: \sigma\subseteq\rho\ \wedge\ \rho
\subseteq\sigma\}$.)\hfill\smallskip

For a Laver-tree $T$, we say $A\subseteq T$ is a {\em{front}\/}
if $\sigma\neq\tau$ in $A$ implies $\sigma\not\subseteq\tau$ and
for all $f\in [T]$ there is an $n\in\omega$ such that $f|_n\in A$.
\hfill\smallskip

The meaning of $p\leq \lval\Phi\rval$ and $p\cap \lval\Phi\rval$
are $U_p\subseteq \lval\Phi\rval$ and $U_p\cap \lval\Phi\rval$,
respectively.\hfill

\begin{enumerate}
\item We say a forcing notion $\Pp$ is {\em{Laver-like}\/} if
there is a $\Pp$-name $\tilde{r}$ for a dominating real such that\\
(i) the complete Boolean
    algebra generated by the family $\{\lval\tilde{r}(i)=n\rval : i,n
    \in\omega\}$ equals r.o.~($\Pp$), and\\
(ii) for each condition $p\in \Pp$ there exists a Laver-tree $T\subseteq
{\omega}^{<\omega}$ so that $$\forall\sigma\in T\Bigl(p(T_{\sigma}):=
\prod\limits_{n\in\omega} \sum\limits_{\tau\in T_{\sigma}}
\{p\cap\lval\tilde{r}|_{\lg (\tau )}=\tau\rval :\lg (\tau )=n\}\in
{\mbox{r.o.~}}(\Pp)\setminus \{\mbox{\bf 0}\}\Bigr).$$
We express this by saying $p(T)\neq\emptyset$ where
$p(T):=p(T_{stem(T)}$).

\item If $\tilde{r}$ is a $\Pp$-name that witnesses that
$\Pp$ is Laver-like, we say that $\Pp$
has {\em{strong fusion}\/} if for countably
many open dense sets $D_n\subseteq \Pp$ and for $p\in \Pp$, there is a
Laver-tree $T$ such that $p(T)\neq\emptyset$ and for each $n$:
$$\{\sigma\in T: p(T)\cap\lval\tilde{r}  |_{\lg (\sigma )}=
\sigma\rval\in D_n\}$$ contains a front.

\item A Laver-like $\Pp$ is {\em{closed under finite changes}\/} if
given a $p\in \Pp$ and Laver trees $T$ and $T'$ so that for all
$\sigma\in T':\ |\mbox{succ}_T(\sigma )\setminus\mbox{succ}_{T'}(\sigma )|
<\omega$, if $p(T)\neq\emptyset$,
then $p(T')\neq\emptyset$, too.
\end{enumerate}

Now we call $\Pp$ a {\em{flexible}\/} forcing notion {\em{iff}\/}
$\Pp$ is Laver-like, has strong fusion and is closed under finite changes.

With this definition we can show (as a further symmetry between
the forcing notions $\M$ and $\dM$),
that the dual-Mathias forcing $\dM$ is flexible.

\begin{lm}
\label{lm:dual-flexible}
The dual-Mathias forcing $\dM$ is flexible.
\end{lm}

\proof By $\dM\approx
{\fT}_{\dM}$ it is enough to prove
that the forcing notion ${\fT}_{\dM}$ is flexible. Let ${\tilde{r}}$
be the canonical ${\fT}_{\dM}$-name for the ${\fT}_{\dM}$-generic
object. By the definition of $\sharp$ and the construction of
${\fT}_{\dM}$, ${\tilde{r}}$ is a name for a dominating real. The
rest of the proof is similar to the proof that Mathias forcing is
%Lemma~3.7
flexible, which is given in \cite{Halb}.
\eop

If all $\Sigma^1_n$-sets in $V$ with parameters in $V\cap W$ have
the Ramsey property ${\cal R}$ or the dual-Ramsey property $\dR$,
we will write $V\models\Sigma^1_n({\cal R})_W$ or
$V\models\Sigma^1_n(\dR)_W$, respectively.
If $V=W$, then we omit the index $W$.
The notations for $\Delta^1_n$-sets and $\Pi^1_n$-sets are
similar. Further ${\cal B}$ stands for the Baire property and
${\cal L}$ stands for Lebesgue measurable.

Now we can prove the following
\begin{thm}
If $V$ is $\Sigma^1_4$-$\dM$-absolute, then $\oV$ is inaccessible in $\LL$.
\end{thm}

\proof To prove the corresponding result for Mathias forcing
(cf.\,\cite{Halb})
% proof of the claim in the proof of Theorem~5.3),
we used only that $\M$ is flexible and that, if $V$ is
$\Sigma^1_4$-$\M$-absolute, then $V\models\Sigma^1_2({\cal
R})$, which is the same as $\Sigma^1_3$-$\M$-absoluteness
(cf.\,\cite{Halb}).
%Theorem~4.1).
Therefore it is enough to prove that
$\Sigma^1_3$-$\dM$-absoluteness implies $\Sigma^1_3$-$\M$-absoluteness.
It follows immediately from Fact~\ref{fct:add-Mathias} that $V\subseteq
V^{\M}\subseteq V^{\dM}$. Now because $\Sigma^1_3$-formulas are upwards
absolute, this completes the proof.
\eop

\section{Iteration of dual-Mathias forcing}

In this section we will build two models in which every
$\Sigma^1_2$-set is dual-Ramsey. In the first model
$2^{\aleph_0}=\aleph_1$ and in the second model
$2^{\aleph_0}=\aleph_2$. With the result that dual-Mathias
forcing has the Laver property we can show that $\Sigma^1_2(\dR)$
implies neither $\Sigma^1_2({\cal L})$ nor $\Sigma^1_2({\cal B})$.

In the sequel we will use the same notations as in
section~\ref{sec:dMathias}.

First we give a result similar to
Theorem~1.15 of \cite{delta1-2}.
\begin{lm}
\label{lm:one-step}
Let $\fD$ be $\fU$-generic over $V$.
If $\fw$ is $\Pad$-generic over $V[\fD]$, then
$V[\fD][\fw]\models\Sigma^1_2(\dR)_V$.
\end{lm}

\proof Let $\Ww$ be the canonical name for the $\Pad$-generic object $\fw$
over $V[{\fD}]$ and let $\varphi (X)$ be a
$\Sigma^1_2$-formula  with parameters
in $V$. By Theorem~\ref{thm:gamefilter}
and Corollary~\ref{cor:pd}, the forcing notion $\Pad$ has pure
decision. So, there exists a $\Pad$-condition $p\in V[\fD]$ with empty
stem (therefore $p\in\fD$), such that {$V[\fD]$}$\models
``p{\force}_{\Pad}\varphi (\Ww)"$ or {$V[\fD]$}$\models
``p{\force}_{\Pad}\neg\varphi (\Ww)"$. Assume the former case holds.
Because $\fw\ceq^* q$ for all
$q\in\fD$, there is an $f\in [\omega]^{<\omega}$ such that
${\fw}\puni\{f\}\ceq p$. By Theorem~\ref{thm:gamefilter} and
Theorem~\ref{thm:proper} we know that if
$X$ is $\Pad$-generic over $V[\fD]$ and $X'\in\partX\cap
V[{\fD}][{\fw}]$, then $X'$ is also $\Pad$-generic over $V[\fD]$. Hence every
$\fw'\ceq {\fw}\puni\{f\}\ceq p$ is $\Pad$-generic over $V[\fD]$ and therefore
$V[\fD][\fw']\models\varphi ({\fw'})$. Because $\Sigma^1_2$-formulas are
absolute, we get $V[\fD][\fw]\models\varphi (\fw')$. So,
$V[{\fD}][{\fw}]\models\exists X\forall Y\in\partX \varphi (Y)$.
The case when {$V[\fD]$}$\models ``p{\force}_{\Pad}\neg\varphi ({\Ww})"$
is similar. Hence, we finally have $V[{\fD}][{\fw}]\models\Sigma^1_2(\dR)_V$.
\eop

\rmk The proof of the analogous result can be found in \cite{delta1-2}.
\hfill\smallskip

Because G\"odel's constructible universe $L$ has a
$\Delta^1_2$-well-ordering of the reals, $L$ is neither a
model for $\Delta^1_2 (\dR)$ nor a model for $\Delta^1_2 ({\cal R})$.
But we can build a model
in which $2^{\aleph_0}=\aleph_1$ and all $\Sigma^1_2$-sets
are dual-Ramsey.

\begin{thm}
\label{thm:iteration}
If we make an $\omega_1$-iteration of dual-Mathias forcing with countable
support starting from $\LL$, we get a model in which every
$\Sigma^1_2$-set of reals is dual-Ramsey and
$2^{\aleph_0}=\aleph_1$.
\end{thm}

\proof Follows immediately from the Fact~\ref{fct:fequi},
Lemma~\ref{lm:one-step}
and the fact that the dual-Mathias forcing is proper.\eop

\rmk The proof of a similar result can be found in \cite{sigma1-2}.
\hfill\smallskip

We can build also a model
in which $2^{\aleph_0}=\aleph_2$ and all $\Sigma^1_2$-sets
are dual-Ramsey.

\begin{thm}
\label{thm:iteration2}
If we make an $\omega_2$-iteration of dual-Mathias forcing with countable
support starting from $\LL$, we get a model in which every
$\Sigma^1_2$-set of reals is dual-Ramsey and
$2^{\aleph_0}=\aleph_2$.
\end{thm}

\proof In \cite{Lori} it is shown that a
$\omega_2$-iteration of dual-Mathias forcing with countable
support starting from $\LL$ yields a model in which
$2^{\aleph_0}=\aleph_2$ and the union of fewer than $\aleph_2$
completely dual-Ramsey sets is completely dual-Ramsey.
Now because each $\Sigma^1_2$-set can be written as the union
of $\aleph_1$ analytic sets (and analytic sets are completely
dual-Ramsey) all the $\Sigma^1_2$-sets are dual-Ramsey.
\eop

\rmk A similar result is true because an $\omega_2$-iteration
of Mathias forcing with countable support starting from $\LL$
yields a model in which $\mathfrak{h}=\aleph_2$ (cf.\,\cite{ShSp}),
and $\mathfrak{h}$ can be considered as
the additivity of the ideal of completely Ramsey null sets
(cf.\,\cite{Plewik}).\hfill\smallskip

For the next result we have to give first the definition of the
Laver property.\hfill\smallskip

A {\em cone\/} $\bar{A}$ is a sequence $\la A_k:k\in\omega\ra$ of
finite subsets of $\omega$ with $|A_k|<2^k$. We say that $\bar{A}$
{\em covers\/} a function $f\in\omega^{\omega}$ if for all $k>0$:
$f(k)\in A_k$. For a function $H\in\omega^{\omega}$, we write
$\Pi H$ for the set $\{f\in\omega^{\omega}:\forall k>0(f(k)<H(k))\}$.
Now a forcing notion $\Pp$ is said to have the {\em Laver property\/}
iff for every $H\in\omega^{\omega}$ in $V$,
$$\one\force_{\Pp} ``\forall f\in\Pi H\,\exists\bar{A}\in V
(\mbox{$\bar{A}$ is a cone covering $f$})".$$

Like Mathias forcing, the dual-Mathias forcing has the Laver property
as well and therefore adds no Cohen reals (cf.\,\cite{tools} and
\cite{book}).

\begin{lm}
\label{lm:Laver}
The forcing notion $\dM$ has the Laver property.
\end{lm}

\proof Given $f,H\in {\omega}^{\omega}$ such
that for all $k>0$: $f(k)<H(k)$,
let $\langle s,X\rangle$ be an $\dM$-condition. Because $\dM$ has pure
decision and $f(1)<H(1)$,
we find a $Y_0\in {\open{s}{X}}$ such that $\langle
s,Y_0\rangle$ decides $f(1)$. Set $s_0:=s$. Suppose we have constructed
$s_n\in\NN$ and $Y_n\in\parto$ such that $s\seg s_n$, $|s_n|=|s|+n$ and
${\open{s_n}{Y_n}}$ is a dual Ellentuck neighborhood.
Choose $Y_{n+1}\in {\open{s_{n}}{Y_n}}$ such that for all $h\in\NN$ with
$s\seg h\ceq s_{n}$ and $\mdom (h)=\mdom (s_{n})$:
$\langle h,Y_{n+1}\rangle$ decides $f(k)$ for all $k<2^{n+1}$.
Further let $s_{n+1}\in\NN$ be such that $s_n\seg s_{n+1}$,
$|s_{n+1}|=|s_n|+1=|s|+n+1$ and $s_{n+1}\seg Y_{n+1}$.
Finally let $Y$ be such that for all $n\in\omega$: $s_n\seg Y$.
Evidently, the $\dM$-condition $\langle s,Y\rangle$ is stronger than
the given $\dM$-condition $\langle s,X\rangle$ (or equal). Now if
$k,n\in\omega$ such that $2^n\leq k<2^{n+1}$, then let $\{h_j:j\leq m\}$
be an enumeration of all $s\seg h\ceq s_n$ with $\mdom (h)=\mdom (s_n)$.
It is clear that $m<2^{2^n}$. Further let $A_k:=\{l\in\omega :
\exists j\leq m(\langle h_j,Y\rangle {\force}_{\dM} f(k)=l)\}$; then
$|A_k|\leq m<2^{2^n}$ and because $2^n\leq k$, we have $|A_k|<2^k$. If
we define $A_0:=\{l\in\omega :\langle s,Y\rangle {\force}_{\dM} f(0)=l\}$
then the $\dM$-condition $\langle s,Y\rangle$ forces that
${\bar{A}}:=\langle A_k: k\in\omega\rangle$ is a cone for $f$.
\eop

Using these results we can prove the following
\begin{thm}
\label{thm:Ramsey-Baire}
$\Sigma^1_2(\dR)$ implies neither $\Sigma^1_2({\cal L})$ nor
$\Sigma^1_2({\cal B})$.
\end{thm}
\proof Because a forcing notion with the Laver property adds no Cohen
reals and because the Laver property
is preserved under countable support
iterations of proper forcings (with the Laver property), in the model
constructed in Theorem~\ref{thm:iteration} no real is Cohen over $\LL$.
Therefore in this model $\Delta^1_2({\cal B})$ fails and because
$\Sigma^1_2({\cal L})$ implies $\Sigma^1_2({\cal B})$
(by \cite{sigma1-2}) also $\Sigma^1_2({\cal L})$
has to be wrong in this model.
\eop

\rmk For the analogous result cf.\,\cite{delta1-2}.

\section{Appendix}

Although the Ramsey property and the dual-Ramsey
property are very similar, we can show that the two Ramsey
properties are different.
\begin{thm}
\label{thm:different}
With the axiom of choice we can construct a set which is Ramsey
but not dual-Ramsey.
\end{thm}
\proof First we construct a set $C\subseteq\reals$ which is not
dual-Ramsey. The relation $``\ass "$ is an
equivalence-relation on $\parto$.
(Remember that $X\ass Y$ if and only if there are $f,g\in
[\omega]^{<\omega}$ such that $X\puni \{f\}\ceq Y$ and $Y\puni \{g\}\ceq
X$.) Now choose from each equivalence
class $X^{*}$ an element $A_X$ and let $h_X:=|f|+|g|$ be of least
cardinality, where $f$ and $g$ are such that $X\puni\{f\}\ceq A_X$ and
$A_X\puni \{g\}\ceq X$.
Further define:
\[ F(X):= \left\{ \begin{array}{ll}
                  1 & \mbox{if } h_X \mbox{ is odd,}\\
                  0 & \mbox{otherwise.}
                  \end{array}\right. \]
Then the set $\{X\in\parto: F(X)=1\}$ is evidently not dual-Ramsey and
therefore also the set $C:=\{x\in\reals :\exists X\in\parto (F(X)=1\wedge
x=\pac (X))\}$ is not dual-Ramsey.\hfill\smallskip

Now define $r:=\{\flat\{k,k+1\}:k\in\omega\}$,
then $\cp (r)=\{\omega\}\not\in\parto$ and hence
$[r]^{\omega}\cap C=\emptyset$. So, the set $C$ is Ramsey.
\eop
We can show that the dual-Ramsey property is stronger than the
Ramsey property.

\begin{lm}
\label{lm:stronger}
If $V\models\Sigma^1_n(\dR)$ then $V\models\Sigma^1_n(\cal R)$.
\end{lm}

\proof Given a $\Sigma^1_n$-formula $\varphi (x)$ with parameters in
$V$. Let $\psi (y)$ be defined as follows.
$$\psi (y)\ \mbox{{\em iff}}\ \ \exists x(x=\mMin (\cp (y))\wedge\varphi
(x)).$$ We see that $\psi (y)$ is also a $\Sigma^1_n$-formula (with the
same parameters as $\varphi$). Now if there is an $X\in\parto$
such that for all $Y\in\partX$,
$\psi (\pac (Y))$ holds, then for all $y\in [x]^{\omega}$ where
$x=\mMin (X)$, $\varphi (y)$ holds. The case where for all $Y\in\partX$,
$\neg\psi (\pac (Y))$ holds, is similar.
\eop

With these results and all the symmetries we found
between the two Ramsey properties and between
the Mathias forcing and the dual-Mathias forcing, it is na\-tural to ask
whether there is a pro\-perty which is equivalent to ``every
$\Sigma^1_2$-set of reals has the dual-Ramsey property". Another
interesting open problem, which surely would give us
a lot of information about the relationship between the two Ramsey
properties, would be the following question:
$$\mbox{Is}\ \Sigma^1_2 ({\cal R})\ \mbox{equivalent to}\ \Sigma^1_2
(\dR)\mbox{?}$$

{\bf Acknowledgements:}\ I would like to thank Pierre Matet
for all the fruitful discussions concerning the results in this paper,
and the referee for valuable comments on an earlier
version of this paper.

\medskip

Eidgen\"ossische Technische Hochschule\\
Departement Mathematik\\
ETH-Zentrum\\
8092 Z\"urich\\
Switzerland\hfill\smallskip

E-mail: halbeis@@math.ethz.ch


\begin{thebibliography}{Maa}

\bibitem[BJ]{book}
{\bo{T.~Bartoszy\'{n}ski and H.~Judah}
    {Set Theory: on the structure of the real line}
    {A.\,K.\,Peters}
    {Wellesley}
    {1995}}

\bibitem[Br]{Brendle}
{\art{J.~Brendle}
     {Combinatorial properties of classical forcing notions}
     {Ann.\,of Pure and Applied Logic}
     {73}
     {1995}
     {143--170}}

\bibitem[CS]{dual}
{\art{T.~J.~Carlson and S.~G.~Simpson}
      {A dual form of Ramsey's Theorem}
      {Adv.\,in Math.}
      {53}
      {1984}
      {265--290}}

\bibitem[El]{Ellentuck}
{\art{E.~Ellentuck}
     {A new proof that analytic sets are Ramsey}
     {J.\,Symb.\, Logic}
     {39}
     {1974}
     {163--165}}

\bibitem[Go]{tools}
{\samp{M.~Goldstern}
     {Tools for your forcing construction}
     {Set Theory of the Reals}
     {H.~Judah}
     {305--360}
     {Israel Mathematical Conference Proceedings}
     {Bar-Ilan University}
     {1993}}

\bibitem[Ha]{Lori}
{\unfart{L.~Halbeisen}
     {On shattering, splitting and reaping partitions}
     {preprint}}

\bibitem[HJ]{Halb}
{\art{L.~Halbeisen and H.~Judah}
     {Mathias absoluteness and the Ramsey property}
     {J.\,Symb.\,Logic}
     {61}
     {1996}
     {177--193}}

\bibitem[HM]{HalMat}
{\prep{L.~Halbeisen and P.~Matet}}

\bibitem[Je 1]{Jech}
{\bo{T.~Jech}
    {Set Theory}
    {Academic Press}
    {London}
    {1978}}

\bibitem[Je 2]{multiple}
{\bo{T.~Jech}
    {Multiple Forcing}
    {Cambridge University Press}
    {Cambridge}
    {1987}}

\bibitem[Ju]{sigma1-2}
{\art{H.~Judah}
     {$\Sigma^1_2$-sets of reals}
     {J.\,Symb.\,Logic}
     {53}
     {1988}
     {636--642}}

\bibitem[JS]{delta1-2}
{\art{H.~Judah and S.~Shelah}
     {$\Delta^1_2$-sets of reals}
     {Ann.\,of Pure and Applied Logic}
     {42}
     {1989}
     {207--223}}

\bibitem[Ka]{Kastanas}
{\art{I.\,G.\,Kastanas}
    {On the Ramsey property for sets of reals}
    {J.\,Symb.\,Logic}
    {48}
    {1983}
    {1035--1045}}

\bibitem[Ke]{Kechris}
{\bo{A.~S.~Kechris}
    {Classical Descriptive Set Theory}
    {Springer-Verlag}
    {New York}
    {1995}}

\bibitem[Ku]{Kunen}
{\bo{K.~Kunen}
    {Set Theory, an Introduction to Independence Proofs}
    {North Holland}
    {Am\-st\-erdam}
    {1983}}

\bibitem[Mt 1]{partitions}
{\art{P.~Matet}
     {Partitions and filters}
     {J.\,Symb.\,Logic}
     {51}
     {1986}
     {12--21}}

\bibitem[Mt 2]{Matet}
{\art{P.~Matet}
     {Happy families and completely Ramsey sets}
     {Arch.\,Math.\,Logic}
     {32}
     {1993}
     {151--171}}

\bibitem[Mt 3]{Matet3}
{\unfart{P.~Matet}
     {Combinatorics and forcing with distributive ideals}
     {to appear in {\em Ann. of Pure and Applied Logic}}}

\bibitem[Ma]{happy}
{\art{A.~R.~D. Mathias}
     {Happy families}
     {Ann.\,Math.\,Logic}
     {12}
     {1977}
     {59--111}}

\bibitem[Pl]{Plewik}
{\art{S.~Plewik}
     {On completely Ramsey sets}
     {Fund.\,Math.}
     {127}
     {1986}
     {127--132}}

\bibitem[SS]{ShSp}
{\unfart{S.~Shelah and O.~Spinas}
     {The distributivity numbers of products
      of ${\cal P}(\omega)$/fin I}
     {preprint}}

\bibitem[Si]{Silver}
{\art{J.~Silver}
     {Every analytic set is Ramsey}
     {J.\,Symb.\,Logic}
     {35}
     {1970}
     {60--64}}

\end{thebibliography}
\end{document}